\documentclass{amsart}
\usepackage{amsmath, amsthm, amssymb}

\theoremstyle{definition}
\newtheorem{defn0}{Definition}[section]
\newtheorem{exa}[defn0]{Example}
\theoremstyle{remark}

\theoremstyle{plain}
\newtheorem{prop0}[defn0]{Proposition}
\newtheorem{thm0}[defn0]{Theorem}

\newtheorem{coro0}[defn0]{Corollary}

\newcommand{\propref}[1]{Proposition~\ref{#1}}
\newcommand{\thmref}[1]{Theorem~\ref{#1}}

\newcommand{\cororef}[1]{Corollary~\ref{#1}}
\newcommand{\exref}[1]{Example~\ref{#1}}

\newcommand{\p}[1]{\mathbb{P}^{#1}}
\newcommand{\codim}[1]{\operatorname{codim} #1}

\begin{document}
\title[Chern numbers of smooth varieties]{Chern numbers of smooth varieties via homotopy continuation and intersection theory}
\author[S. Di Rocco]{Sandra Di Rocco}
    \address{Department of Mathematics,
             KTH, 100 44 Stockholm, Sweden}
    \email{dirocco@math.kth.se}
    \urladdr{http://www.math.kth.se/$\sim$sandra}
\author[D. Eklund]{David Eklund}
    \address{Department of Mathematics,
               KTH, 100 44 Stockholm, Sweden}
    \email{daek@math.kth.se}
     \urladdr{http://www.math.kth.se/$\sim$daek}
\author[C. Peterson]{Chris Peterson}
    \address{Department of
             Mathematics, Colorado State University, Fort Collins, CO 80523}
    \email{peterson@math.colostate.edu}
    \urladdr{http://www.math.colostate.edu/$\sim$peterson}
\author[A.J. Sommese]{Andrew  J. Sommese}
    \address{Department of Mathematics,
             University of Notre Dame, Notre Dame, IN 46556}
    \email{sommese@nd.edu}
    \urladdr{http://www.nd.edu/$\sim$sommese}
\thanks{The first  author was partially supported by V.R. grant NT:2006-3539. The second author was mainly supported by
Knut och Alice Wallenbergs Stiftelse and partially supported by CIAM, KTH. The third author
was partially supported by NSF  MSPA-MCS-0434351 and
AFOSR-FA9550-08-1-0166. The fourth author was partially
supported by the Duncan Chair of the University of Notre Dame, NSF
DMS-0410047 and NSF DMS-0712910.}

\begin{abstract}
Homotopy continuation provides a numerical tool for computing the equivalence of a smooth variety in an intersection product. Intersection theory provides a theoretical tool for relating the equivalence of a smooth variety in an intersection product to the degrees of the Chern classes of the variety. A combination of these tools leads to a numerical method for computing the degrees of Chern classes of smooth projective varieties in $\p{n}$. We illustrate the approach through several worked examples.
\end{abstract}
\maketitle
\section{Introduction}
In this paper we describe a numerical method for computing the degrees of the Chern classes of a smooth projective variety. The main tools that will be used are drawn from intersection theory and numerical algebraic geometry. In particular, we will be using a refined B\'{e}zout's theorem \cite{FOV,Fu,Vo} together with ideas related to numerical polynomial algebra such as homotopy continuation, monodromy, and the numerical decomposition of zero sets \cite{Li, SVW, SW, St}. In addition, we will need to numerically determine the degree of a certain residual zero-scheme \cite{BPS, DZ, Le, Ley, LVZ}.  Implementations of many of the computational tools that will be used can be found in freely available software packages (e.g. \cite{BHSW, Co, GS, GFS, GKKTFM, LLT, SJ, Ve, WSMMW, Ze}).  The algorithms developed in this paper primarily use the software package {\it Bertini} with help from the software package {\it Macaulay 2} \cite{BHSW,GS}.

From the generators of the ideal of a smooth projective variety, Chern numbers can be computed through a purely symbolic computation. This naturally leads one to ask if a numerical approach is useful or needed. A partial answer is that in the numerical approach described in this paper, the algorithms work equally well whether the generators are sparse or dense and whether they have rational, algebraic or transcendental coefficients. This is not the case with purely symbolic methods.  Furthermore, a surprising amount of information can be extracted from an ideal even in situations where the generators have inaccuracies in their coefficients \cite{HMPS, HSS, NSW}. Another important feature, that will likely play a dominant role in the future, is that homotopy continuation algorithms parallelize well. This will allow such algorithms to take advantage of the trend of more processors on a chip. These features lead us to believe that the introduction of a numerical approach is indeed useful and will allow experimentation in settings outside the domain of purely symbolic methods.

\subsubsection*{Acknowledgments:} The authors would like to thank Frank-Olaf Schreyer and Charles Wampler for stimulating conversations and observations.
\section{intersection Theory and Homotopy Continuation}
In this section, we describe the main tools that will be used from intersection theory and from numerical algebraic geometry. Let $X_1,\dotsc,X_r$ be a set of hypersurfaces in $\p{r}$. Let $Z$ be a smooth connected component of  the scheme defined by $X_1 \cap \dotsc \cap X_r$. We present a well known theorem from intersection theory that relates the equivalence of $Z$ in the intersection product $X_1 \cdot \dotsc \cdot X_r$ to a formula involving the total Chern classes of certain bundles related to $X_1,\dotsc,X_r, Z$, and $\p{r}$. With mild assumptions on the scheme determined by the hypersurfaces, this theorem can be translated into a numerical condition involving the degrees of the Chern classes of $Z$, the degrees of the $X_i$ and the degree of a certain residual zero-scheme. Tools from numerical algebraic geometry can be used to compute the degree of the residual zero-scheme. By varying the degrees of the $X_i$, we obtain multiple numerical conditions that can be used to uniquely determine the degrees of the Chern classes of $Z$.

\subsection{Intersection Theory}
Intersection theory has been primarily developed from the two different viewpoints found in \cite{FOV, Fu, Vo}. Though each viewpoint seems to have its own advantages and disadvantages, Leendert van Gastel showed that the two are closely related \cite{vG}. While either version would be adequate for our purposes, we follow the approach of Fulton-MacPherson as described in \cite{Fu}. We first need to fix the notation and definitions that will be used in this section.
\begin{defn0} Let $Z$ and $X$ be projective varieties in $\mathbb P^r=\mathbb P_{\mathbb C}^r$ with $Z\subseteq X$.
\begin{enumerate}
\item Let $A_*(X)= \bigoplus_k A_k(X)$ where $A_k(X)$ denotes the Chow group of $k$-dimensional cycles on $X$ modulo rational equivalence.\item For $\alpha \in A_{\ast}(X)$, ${\{\alpha \}}_0$  denotes the component of $\alpha$ in $A_0(X)$. Let $\alpha \in A_k(X)$, $\alpha=\sum_{i=1}^m a_iV_i$ where $a_1,\dots,a_m \in \mathbb{Z}$ and $V_1,\dots,V_m \subseteq X$ are irreducible subvarieties. The degree of $\alpha$ is defined by $\deg{\alpha}=\sum_i a_i \deg{V_i}$.
\item $N_{X} \p{r}$ denotes the normal bundle of $X$ in $\mathbb P^r$, $T_X$ denotes the tangent bundle of $X$ and $T_X |_Z$ denotes its restriction to $Z$.
\item Given a vector bundle $E$ of rank $n$ on a smooth variety, $c_i(E)$ denotes the $i^{th}$ Chern class of $E$ and $c(E)=1+c_1(E)+\dots + c_n(E)$  denotes the total Chern class of $E$.
\item The $i^{th}$ Chern class of a smooth variety $X$ is defined as $c_i(T_X)$ and $\deg(c_i(T_X))$ is called a Chern number of $X$.
\end{enumerate}
\end{defn0}

Let $X_1,\dotsc,X_r$ be hypersurfaces in $\p{r}$ and let $Z$ be a smooth connected component of  $X_1 \cap \dotsc \cap X_r$. In \cite{Fu} Fulton considers the intersection product $X_1 \cdot \dotsc \cdot X_r \in A_0(X_1 \cap \dotsc \cap X_r)$ and
defines $(X_1 \cdot \dotsc \cdot X_r)^Z \in A_0(Z)$, the \emph{equivalence} of $Z$ for $X_1 \cdot \dotsc \cdot X_r$, as the part of the intersection product supported on $Z$. The following is a simplified version of Proposition 9.1.1  of \cite{Fu} specialized to our setting.

\begin{prop0} \emph{(Fulton)} \label{prop:equivalence}
Let $X_1,\dotsc,X_r$ be hypersurfaces in $\p{r}$ and let $Z$ be a smooth connected component of $X_1 \cap \dotsc \cap X_r$. Let $N_i$ be the restriction of $N_{X_i} \p{r}$ to $Z$. Then
\begin{equation}\label{formula}(X_1 \cdot \dotsc \cdot X_r)^Z = {\{(\Pi_{i=1}^r c(N_i))c(T_{\p{r}} |_Z)^{-1}c(T_Z)\}}_0.\end{equation}
\end{prop0}

We introduce notation for the purpose of stating a useful corollary of \propref{prop:equivalence}. Let $n_i=\deg{X_i}$, let $\dim{Z}=n$,  and let $c_0,\dots,c_n$ be the Chern classes of $Z$. The $k^{th}$ elementary symmetric function in $n_1,\dots, n_r$ will be denoted by $\sigma_k$. In other words,
$$\sigma_0=1, \quad \sigma_1 = \sum_{i} n_i, \quad \sigma_2=\sum_{i<j} n_i n_j, \quad \sigma_3=\sum_{i<j<k} n_i n_j n_k, \quad \dotsc$$ If we let
$$a_i=\sum_{j=0}^{n-i} (-1)^j \binom{r+j}{j} \sigma_{n-i-j}$$
then we obtain the following identity as a direct corollary of \propref{prop:equivalence}:
\begin{coro0} \label{coro:formula}
$\deg{(X_1 \cdot \dotsc \cdot X_r)^Z} = \sum_{i=0}^n a_i \deg{c_i}.$
\end{coro0}
\begin{proof}
Let $H$ be the hyperplane class in $\p{r}$ and let $H|_Z$ be its restriction to $Z$. Because  $c(T_{\p{r}} |_Z)=(1+H|_Z)^{r+1}$, $$c(T_{\p{r}} |_Z)^{-1}=\sum_{j=0}^n(-1)^j {r+j\choose j}H|_Z^j.$$
As cycle classes on $\p{r}$, $N_{X_i} \p{r}=X_i=n_iH$. Hence $c(N_i)=1+n_iH|_Z$ and we get that $\Pi_{i=1}^r (1+c(N_i))=\sum_{k=0}^n \sigma_k H|_Z^k$. Thus the formula (\ref{formula}) can be written as
$$(X_1 \cdot \dotsc \cdot X_r)^Z = {\left\{(\sum_{k=0}^n \sigma_k H|_Z^k)(\sum_{j=0}^n(-1)^j {r+j\choose j}H|_Z^j)(\sum_{i=1}^n c_i) \right\}}_0.$$
The statement now follows since $\deg{c_i}=\deg{(c_i H|_Z^{n-i})}$ and the component of $(X_1 \cdot \dotsc \cdot X_r)^Z$ in $A_0(Z)$
is given by the sum of terms where $k=n-i-j$.
\end{proof}
The classical version of B\'{e}zout's theorem (Proposition 8.4 of \cite{Fu}) states that $$\deg{(X_1 \cdot \dotsc \cdot X_r)}=\prod_{i=1}^r n_i.$$ A refined version of B\'{e}zout's theorem (Proposition 9.1.2 of \cite{Fu}) ties these results together with an extended B\'{e}zout formula.
\begin{prop0} \label{eq:sum} Suppose $X_1 \cap \dotsc \cap X_r$ consists of a connected component $Z$ and a finite set $S$. For $p \in S$ let  $m_p=i(p,X_1 \cdot \dotsc \cdot X_r;\p{r})$ denote the intersection multiplicity of $p$ in $X_1 \cdot \dotsc \cdot X_r$. Then
$\deg{(X_1 \cdot \dotsc \cdot X_r)^Z}+\sum_{p \in S}m_p=\prod_{i=1}^r n_i.$
\end{prop0}
If $n_1,\dotsc,n_r$ and $\sum_{p \in S} m_p$ are known, then we can use \propref{eq:sum} to solve for  $\deg{(X_1 \cdot \dotsc \cdot X_r)^Z}$. If $Z$ is smooth, \cororef{coro:formula} provides a linear relation among the degrees of the Chern classes of $Z$. We summarize these observations in the following theorem implicit in Fulton \cite{Fu}:
\begin{thm0} \label{num:criterion}
 Let $F_1, \dotsc, F_r$ be homogeneous forms corresponding to hypersurfaces $X_1, \dotsc, X_r \subset \p{r}$. Let $n_i=\deg{X_i}$, let $\sigma_k$ be the $k^{th}$ elementary symmetric polynomial in $n_1,\dots, n_r$ and let $\vec{\mathbf A}=[a_0,\dotsc a_n]$ with $a_k=\sum_{i=0}^{n-k} (-1)^i \binom{r+i}{i} \sigma_{n-k-i}$. Let $Z$ be a smooth connected $n$-dimensional scheme with Chern classes $c_0,\dotsc, c_n$ and let $\vec{\mathbf C}=[\deg{c_0}, \dotsc, \deg{c_n}]$. If the subscheme of $\p{r}$ defined by the ideal $(F_1, \dotsc, F_r)$ is a disjoint union of $Z$ and a (possibly empty) zero-scheme $S$, then $\vec{\mathbf A} \cdot \vec{\mathbf C}= \sigma_r -\deg{S}.$
\end{thm0}
In order to compute the degrees of the Chern classes of an $n$-dimensional smooth variety, it is enough to determine $n+1$ independent linear relations that they satisfy. \thmref{num:criterion} provides a mechanism for producing the linear relations provided we are able to find a sufficient number of $r$-tuples of homogeneous forms that satisfy the conditions of \thmref{num:criterion} and provided we have a method for computing $\deg{S}$. In the next two sections, we present a collection of tools that allow us to complete the latter task.

\subsection{Homotopy Continuation}
In homotopy continuation, a polynomial ideal, $I$,
is cast as a member of
a parameterized family of polynomial ideals one of which has known
isolated solutions. Each of the known isolated solutions is tracked through a
predictor/corrector method to a point which lies numerically close to the algebraic set $V(I)$ determined by $I$. These
points can then be refined to lie within a prescribed tolerance of $V(I)$. An introduction to general
continuation methods can be found in \cite{AG}.
Through the basic algorithms of numerical algebraic geometry \cite{SW},
 from $I$  it is then possible to
produce a collection of subsets of points such
that:
\begin{itemize}
\item The subsets are in one to one correspondence with the
irreducible components of the algebraic set $V(I)$.
\item The points in a subset all lie within a prescribed tolerance of
the irreducible component to which it corresponds.
\item The number of
points in the subset is the same as the degree of the irreducible
component.
\item The subset is a numerical approximation of the
intersection of the irreducible component with a known linear space
of complementary dimension.
\end{itemize}
Note that the decomposition above allows one to determine both the dimension and degree of each algebraic variety appearing in the decomposition of an algebraic set.

\subsection{Degree of a Zero Dimensional Scheme}
To each
irreducible component of an algebraic set, one can use the
defining set of polynomials to attach a positive integer, called the {\it multiplicity}, that
determines roughly how many times the component should be counted
in a computation.
In \cite{DZ}, Dayton and Zeng study the multiplicity inspired, to a large degree, by Macaulay's inverse
systems approach.  They provide an algorithm which yields as
output the multiplicity of isolated solutions. This is essentially done by
counting how many partial derivatives of the polynomials are forced
to be zero. In \cite{BPS}, an alternate approach is presented,
inspired, to a large degree, by
certain Gr\"obner basis calculations coupled with a fundamental result
of Bayer and Stillman on regularity \cite{BS}. Both of these algorithms have been implemented in \emph{Bertini} \cite{BHSW}. Thus, there is an implemented algorithm in place that allows one to determine the degree of a zero-dimensional scheme.

\section{Algorithms}
In this section, we present the pseudocode for  three algorithms. The first algorithm computes the equivalence of a connected scheme $Z$ in an intersection product. The second algorithm determines a linear relation satisfied by the degrees of the Chern classes of a smooth connected scheme. The third algorithm computes the degrees of the Chern classes of a smooth connected scheme. For each of the algorithms, the input is a set of non-zero homogeneous generators for an ideal $I$. It is an assumption of the algorithms that the scheme determined by $I$ is a disjoint union of $Z$ and a zero-scheme $S$. The second and third algorithms have the additional assumption that $Z$ is smooth.
\\ \\* \\*
{\bf Algorithm~1.} {\it Equivalence\_of\_Z} $(\{F_1,F_2,\dots ,F_r\}; D)$ \\ \\*
\underline{\bf Input}: A set of $r$ homogeneous polynomials $\{F_1,F_2,\dots, F_r\}\subset \mathbb C[z_0, z_1, \dots ,z_r]$. The polynomials should generate an ideal $I$ whose corresponding scheme is the disjoint union of a connected scheme $Z$ and a possibly empty zero-scheme $S$. \\ \\*
\underline{\bf Output}: $D=\deg{(X_1 \cdot \dotsc \cdot X_r)^Z}$ where $X_i$ is the hypersurface corresponding to $F_i$.\\ \\*
\underline{\bf Algorithm}:\\
\hspace*{0.15in}Determine the support of $S$.\\
\hspace*{0.15in}Determine the multiplicity of each point in the support of $S$.\\
\hspace*{0.15in}Add up the multiplicities of the points in the support of $S$ and store in $\mu_S$.\\
\hspace*{0.15in}For each $i$, determine the degree of $F_i$ and store in $n_i$.\\
\hspace*{0.15in}Compute $T=\prod_i n_i$.\\
\hspace*{0.15in}Compute $T-\mu_S$ and store in $D$.
\\ \\* \\*
{\bf Algorithm~2.} {\it Linear\_Relation\_On\_Chern\_Numbers} $(\{F_1,F_2,\dots ,F_r\}; \vec{\mathbf A}, D)$ \\ \\*
\underline{\bf Input}: A set of $r$ homogeneous polynomials $\{F_1,F_2,\dots F_r\}\subset \mathbb C[z_0, z_1, \dots ,z_r]$. The polynomials should generate an ideal $I$ whose corresponding scheme is the disjoint union of $Z$ and $S$, where $Z$ is a smooth connected scheme and $S$ is a possibly empty zero-scheme.\\ \\*
\underline{\bf Output}: $\vec{\mathbf A}=[a_0,\dots, a_n]$ and $D\in \mathbb Z,$ where $n$ denotes the dimension of $Z$. If $c_i$ denotes the $i^{th}$ Chern class of $Z$ and $\vec{\mathbf C}=[\deg{c_0},\dots, \deg{c_n}]$ then the linear relation is $\vec{\mathbf A}\cdot\vec{\mathbf C}=D$.\\ \\*
\underline{\bf Algorithm}:\\
\hspace*{0.15in}Determine the dimension of $Z$ and store in $n$.\\
\hspace*{0.15in}Compute {\it Equivalence\_of\_Z} $(\{F_1,F_2,\dots ,F_r\}; D)$.\\
\hspace*{0.15in}Compute the elementary symmetric functions $\sigma_0,\sigma_1, \dots, \sigma_n$ of $n_1, \dots, n_r$.\\
\hspace*{0.15in}Compute $\vec{\mathbf A}=[a_0, a_1, \dots, a_n]$ where $a_k=\sum_{i=0}^{n-k} (-1)^i \binom{r+i}{i} \sigma_{n-k-i}$.
\\ \\* \\*
{\bf Algorithm~3.} {\it Chern\_Numbers} $(\{F_1,F_2,\dots ,F_t\}; \vec{\mathbf C})$\\ \\*
\underline{\bf Input}: A set of $t$ homogeneous polynomials $\{F_1,F_2,\dots F_t\}\subset \mathbb C[z_0, z_1, \dots ,z_r]$. The polynomials should generate an ideal $I$ whose corresponding scheme is the disjoint union of $Z$ and $S$, where $Z$ is a smooth connected scheme and $S$ is a possibly empty zero-scheme.
\\ \\*
\underline{\bf Output}: $\vec{\mathbf C}=[\deg{c_0},\dots, \deg{c_n}]$ where $n$ denotes the dimension of $Z$ and $c_i$ denotes the $i^{th}$ Chern class of $Z$.\\ \\*
\underline{\bf Algorithm}:\\
\hspace*{0.15in}Determine the dimension of $Z$ and store in $n$. \\
\hspace*{0.15in}Determine the degrees of $F_1,\dots,F_t$, and store the maximal degree as $b$.\\
\hspace*{0.15in}Set $l_c=b$ for $c=1,\dots, n+1$.\\
\hspace*{0.15in}For $i=1$ to $n+1$\\
\hspace*{0.24in}Choose $r$ random elements from $I$ of degrees $l_1, l_2, \dots, l_r$. \\
\hspace*{0.24in}Store the random elements as $G_1, G_2, \dots, G_r$.\\
\hspace*{0.24in}Compute {\it Linear\_Relation\_On\_Chern\_Numbers} $(\{G_1,G_2,\dots ,G_r\}; \vec{\mathbf A}_i, D_i)$.\\
\hspace*{0.24in}Let $l_i=l_i+1$.\\
\hspace*{0.15in}Next $i$ \\
\hspace*{0.15in}Build the $(n+1)\times (n+1)$ matrix $M$ whose $i^{th}$ row is $\vec{\mathbf A}_i$.\\
\hspace*{0.15in}Build $\vec{\mathbf D}=[D_1,\dots, D_{n+1}]$.\\
\hspace*{0.15in}Set up the linear system $M\vec{\mathbf C}=\vec{\mathbf D}$, with $\vec{\mathbf C}$ and $\vec{\mathbf D}$ as column vectors.\\
\hspace*{0.15in}Solve the linear system for $\vec{\mathbf C}$. 
\\ \\*

There are two potential problems that should be addressed in the use of Algorithm~3. The first potential problem is whether the matrix $M$ has full rank. The second potential problem is whether the random $r$-tuples $G_1,G_2, \dots, G_r$, produced in such a simple manner from $I$, satisfy the input requirements of Algorithm~2. These two questions are answered in \propref{prop:det} and \cororef{coro:Bertinimaxdeg}, respectively.
\begin{prop0} \label{prop:det}
The matrix $M$ defined in Algorithm~3 satisfies $\det(M) = \pm 1$.
\end{prop0}
\begin{proof}
View $M$ as a function of the maximal degree $b$. Let $\sigma_i$ be the $i^{th}$ elementary symmetric polynomial in $r$ variables $b_1,\dots,b_r$ and let  $$\sigma_i^j=\sigma_i(b_1+1,b_2+1,\dots,b_j+1,b_{j+1},\dots,b_r).$$ Define $\tau_i^j \in \mathbb{Z}[b]$ by $\tau_i^j=\sigma_i^j(b,b,\dots,b)$, $\tau_i^0=\sigma_i(b,b,\dots,b)$. Define a $(n+1) \times (n+1)$ matrix $N(b)$ by
\begin{displaymath}
N=\left(
\begin{array}{ccccc}
\tau_n^0 & \tau_{n-1}^0 & \dots & \tau_1^0 & 1 \\
\tau_n^1 & \tau_{n-1}^1 & \dots & \tau_1^1 & 1 \\
& & \vdots \\
\tau_n^n & \tau_{n-1}^n & \dots & \tau_1^n & 1 \\
\end{array}
\right).
\end{displaymath}
Observe first that $M$ can be factored as $M=NR$, where
\begin{displaymath}
R=\left(
\begin{array}{ccccccc}
1 & 0 & 0 & \dotsc & 0 & 0 \\
-(r+1) & 1 & 0 & \dotsc & 0 & 0 \\
\binom{r+2}{2} & -(r+1) & 1 & \dotsc & 0 & 0\\
& & \vdots \\
(-1)^{n} \binom{r+n}{n} & (-1)^{n-1} \binom{r+n-1}{n-1} & \dotsc & \binom{r+2}{2} & -(r+1) & 1
\end{array}
\right).
\end{displaymath}
Since $\det(R)=1$, we have $\det(M)=\det(N)$.

One checks that $\det (N(0))=\pm 1$. We shall show that the determinant of $N(b)$ does not depend on $b$ by showing that $$\frac{d(\det (N))}{db}=0.$$ We shall use the following
fact:
$$\frac{d\tau_i^j}{db}=(r-i+1)\tau_{i-1}^j.$$ Let $C_i$ be the $i^{th}$ column of $N$, $N=(C_0,C_1,\dots,C_n)$. Then we have that $$\frac{d(\det (N))}{db}= \sum_
{i=0}^n \det (C_0,C_1,\dots,C_{i-1},\frac{dC_i}{db},C_{i+1},\dots,C_n),$$ but each term in the above sum is zero since $\frac{dC_i}{db}=(r-n+i+1)C_{i+1}$ for $0 \leq i \leq n-1$ and $\frac{dC_n}{db}=0$.
\end{proof}

Let $I$ be an ideal defining a scheme $X \subset \p{r}$ and let $I_X$ be the homogeneous ideal of $X$. For $d \in \mathbb{N}$, we say that $X$ is cut out scheme theoretically in degree $d$ if the saturation of the ideal generated by $I(d)$ is equal to $I_X$.

\begin{prop0} \label{prop:Bertini}
Let $I$ be an ideal defining a scheme in $\p{r}$ that is the disjoint union of a smooth connected scheme $Z$ and a possibly empty zero-scheme $S$. Let $G_1,\dots,G_k \in I$, $1 \leq k \leq r$, generate an ideal that defines a scheme $X$, such that the singular locus $X_{\emph{sing}}$ satisfies $\codim{(X_{\emph{sing}}\setminus S)} > k$, and $\codim{(X \setminus Z)} \geq k$. Let $d \in \mathbb{N}$ be such that $Z \cup S$ is cut out scheme theoretically in degree $d$ and fix integers $n_{k+1},\dots,n_r$ with $n_i \geq d$ for all $i$. If $G_{k+1},\dots,G_r$ are general forms in $I$ with $\deg{G_i}=n_i$, then the ideal $J=(G_1, G_2, \dots, G_r)$ defines the disjoint union of $Z$ and a possibly empty zero-scheme $S'$. If $S$ is non-singular or empty, then $S'$ is non-singular.
\end{prop0}
\begin{proof}
Suppose $k \neq r$. A generic form $G_{k+1} \in I$ of degree $n_{k+1}$ cuts down the dimension of every component of $X \setminus (Z \cup S)$ as well as every component of $X_{\emph{sing}} \setminus S$. By induction, the scheme $Y$ defined by $J$ is a disjoint union of $Z$ and a zero-scheme $S'$, and $Y$ is nonsingular away from $S$. In particular $Y$ is non-singular on $Z$.
\end{proof}

\begin{coro0} \label{coro:Bertiniregular}
Suppose $I$ and $d$ are as in \propref{prop:Bertini}. If $G_1,\dots,G_k \in I$ is a regular sequence in $I$ such that the scheme $X$ defined by the ideal $(G_1,\dots,G_k)$ is generically smooth, then for general forms $G_{k+1},\dots,G_r$ with $\deg{G_i} =n_i \geq d$, the ideal $J=(G_1, G_2, \dots, G_r)$ defines the disjoint union of $Z$ and a possibly empty zero-scheme $S'$.
\end{coro0}
\begin{proof}
Every component of $X$ has codimension $k$ and is of multiplicity 1. We may therefore apply \propref{prop:Bertini}.
\end{proof}

The following corollary implies that for generic choices of the
forms $G_1,\dots,G_r$ in Algorithm~3, the assumptions on the input of Algorithm~2 are satisfied.
\begin{coro0} \label{coro:Bertinimaxdeg}
Let $I\neq 0$ be as in \propref{prop:Bertini} and let $\{F_1,\dots,F_t\}$ be a set of non-zero generators of $I$. Put $d=\max \{\deg{F_1},\dots,\deg{F_t}\}$ and fix integers $n_1,\dots,n_r$ with $n_i \geq d$ for all $i$. If $G_1,\dots,G_r$ are general forms in $I$ with $\deg{G_i}=n_i$, then the ideal $J=(G_1, G_2, \dots, G_r)$ defines the disjoint union of $Z$ and a possibly empty zero-scheme $S'$.
\end{coro0}
\begin{proof}
Observe that the ideal generated by $I(d)$ is equal to $\oplus_{e \geq d} I(e)$. Since the ideal $\oplus_{e \geq d} I(e)$ and $I$ have the same saturation, namely the homogeneous ideal of $Z \cup S$, the scheme $Z \cup S$ is cut out scheme theoretically in degree $d$. 
\end{proof}

For reasons of efficiency it is desirable to keep the degrees of the forms $G_1,\dots,G_r$ in Algorithm~3 as low as possible. We give two examples that have bearing on the question of how much
\cororef{coro:Bertinimaxdeg} can be strengthened.

\begin{exa} \label{ex:conic}
Consider a conic $Z \subset \p{3}$. The conic is cut out by a plane $P$ and a quadric $Q$. The corresponding forms of degree 1 and 2 generate the homogeneous ideal $I$ of $Z$. In this case we cannot choose the generic forms $G_1,G_2,G_3 \in I$ of degrees $(3,1,1)$. A cubic and two planes, all three containing the conic, will intersect in the conic union a line, since the planes are necessarily both equal to $P$. The residue in this case is thus a line and not a finite scheme as required.
\end{exa}

\begin{exa}
This example is along the same lines as \exref{ex:conic}.
Let $A$ be a $2 \times 3$ matrix of general linear forms in $\mathbb{C}[x_0,\dots,x_4]$. Let $F_2,F_3,F_4$ be the $2 \times 2$ minors of $A$ and let $F_1$ be a general form of degree 3. The ideal $I=(F_1,F_2,F_3,F_4)$ is the homogeneous ideal of a curve $Z$ in $\p{4}$ and the minors $F_2,F_3,F_4$ define a surface. We thus have a minimal generating set of $I$ of degrees $(3,2,2,2)$ but 4 general forms of degrees $(4,2,2,2)$ will define a union of the curve $Z$ and another curve, violating the input requirements of Algorithm~2.
\end{exa}

\section{Examples}
In this section we present several examples of computations of Chern numbers which illustrate the algorithms of the previous section.

\subsection{Curves} We consider the case of a smooth connected curve in $\p{r}$ (\cite{Fu} Example 9.1.1).
If $Z$ is a smooth curve of genus $g$ in $\p{r}$ then the first Chern class of $Z$ is $-K_{Z}$ and $\deg{(-K_Z)}=2-2g$. If $(F_1,F_2,
\dots, F_r)\subset \mathbb C[z_0, z_1, \dots ,z_r]$ is a homogeneous ideal $I$ whose corresponding scheme is the disjoint union of the curve $Z$ and a zero-scheme $S$ then \cororef{coro:formula} leads to
$$\deg{(X_1 \cdot \dotsc \cdot X_r)^Z}=(n_1+ \dotsc + n_r -(r+1))\deg{Z} +2-2g$$
where $X_i$ denotes the hypersurface corresponding to $F_i$ and $n_i=\deg{X_i}$.
From \thmref{num:criterion}, we can also write the equation as
\begin{equation} \label{eq:curves2}
\prod_in_i - \deg{S}=(n_1+ \dotsc + n_r -(r+1))\deg{Z} +2-2g.
\end{equation}

\begin{exa}
The homogeneous ideal of the twisted cubic curve in $\p{3}$ is $I=(x^2-wy, y^2-xz, wz-xy) \subset \mathbb C[w,x,y,z]$. It is well known that this curve has degree 3 and genus 0. If we choose $F_1, F_2,F_3\in I$ of degrees $(2,2,2)$ then the numerical irreducible decomposition implemented in \cite{BHSW} determines that the corresponding scheme consists of a degree 3 curve and no additional points ($S$ is empty). Equation (\ref{eq:curves2}) gives the relation $$2\cdot2\cdot2-0=(2+2+2-(3+1))\cdot 3 + 2-2g$$ which we can solve to get $g=0$.

If we choose $F_1, F_2,F_3\in I$ of degrees $(2,2,3)$ then we obtain a degree 3 curve and one additional point of multiplicity 1. Equation (\ref{eq:curves2}) gives the relation $$2\cdot2\cdot3-1=(2+2+3-(3+1))\cdot 3 + 2-2g$$ again leading to $g=0$.

If we did not know the degree of $Z$, the two computations in this example would yield $$2\cdot2\cdot2-0=(2+2+2-(3+1))\deg{Z} + 2-2g,$$ $$2\cdot2\cdot3-1=(2+2+3-(3+1))\deg{Z} + 2-2g.$$ The unique solution to these two equations is $\deg{Z}=3$ and $g=0$.
\end{exa}

\subsection{Surfaces} Now consider the case where $Z$ is a smooth connected surface  in $\p{r}$ (\cite{Fu} Example 9.1.5). Let $I=(F_1,F_2,
\dots, F_r)\subset \mathbb C[z_0, z_1, \dots ,z_r]$ be a homogeneous ideal. If the scheme corresponding to $I$ is the disjoint union of the surface $Z$ and a zero-scheme $S$ then \cororef{coro:formula} and \thmref{num:criterion} lead to the equation
\begin{equation} \label{eq:surfaces}
\prod_in_i - \deg{S}=\deg{(X_1 \cdot \dotsc \cdot X_r)^Z}=a_2\deg{c_2}+a_1\deg{c_1}+ a_0 \deg{Z},
\end{equation}
where
$$a_2=1,\ \ a_1=\sum_{i=1}^r n_i -(r+1) \ \ \textrm{and} \ \ a_0=\sum_{i<j} n_i n_j -(r+1)\sum_{i=1}^r n_i+\binom{r+2}{2}.$$
\begin{exa} \label{ex:surface}
In characteristic zero and up to standard modifications, the Horrocks-Mumford bundle $E$ is the only known indecomposable rank 2 vector bundle on $\p{4}$ \cite{HM}.
For a general section $s$ of $E$, the zero set $V(s)$ is a smooth surface in $\p{4}$ called a Horrocks-Mumford surface. We shall compute the degree of a smooth Horrocks-Mumford surface $Z$ as well as the degrees of its Chern classes $c_1$ and $c_2$. It is known that $Z$ is an Abelian surface of degree 10 \cite{HM}. Hence $\deg{c_0}=10$ and $c_1=c_2=0$. The homogeneous ideal $I$ of $Z$ is generated by three quintics and fifteen sextics and $Z$ is cut out scheme theoretically by three quintics and one sextic. Generators for $I$ can be found by a variety of methods, e.g. by a Beilinson monad \cite{DES}. In this example we do not follow Algorithm~3 in detail but the example is covered by \propref{prop:Bertini}.

The following table shows the number of solutions resulting from zero-dimensional, 128-bit precision runs in \cite{BHSW} on 4 random elements of $I$ for various choices of degrees. We let $(n_1,n_2,n_3,n_4)$ denote the degrees of these elements.
\begin{displaymath}
\begin{tabular}{ccccc}
$(n_1,n_2,n_3,n_4)$ & \vline & non-singular points & \vline & singular points \\ \hline
$(5,5,5,6)$ & \vline & 0 & \vline & 750 \\
$(5,5,6,6)$ & \vline & 40 & \vline & 860 \\
$(5,6,6,6)$ & \vline & 100 & \vline & 980 \\
\end{tabular}
\end{displaymath}
In each of these computations the zero set consists of the Abelian surface together with a finite set $S$ of points of multiplicity one. In each row of the table the entry in the middle column is the cardinality of $S$ and the entry in the rightmost column is the B\'{e}zout number $\prod_{i=1}^4 n_i$ minus the entry in the middle column, i.e. $$n_1n_2n_3n_4-\sum_{p \in S} m_p.$$ By \propref{eq:sum}, this is also the degree of the equivalence of $Z$. Evaluating (\ref{eq:surfaces}) in each case gives:
\begin{displaymath}
\begin{array}{l}
\deg{c_2}+16 \deg{c_1}+ 75 \deg{Z}=750 \\
\deg{c_2}+17 \deg{c_1}+ 86 \deg{Z}=860 \\
\deg{c_2}+18 \deg{c_1}+ 98 \deg{Z}=980.
\end{array}
\end{displaymath}
This system has the unique solution $\deg{Z}=10$, $\deg{c_1}=\deg{c_2}=0$.
\end{exa}
\subsection{Higher Dimensional Varieties} If $Z$ is a smooth $n$-dimensional variety  in $\p{r}$ then the algorithm proceeds in a similar manner to the curve and surface cases. We use homotopy continuation to determine the equivalence of $Z$ in various intersection products. This combines with the formulas appearing in \cororef{coro:formula} and \thmref{num:criterion} to produce a linear system involving the degrees of the Chern classes of $Z$. Finally, we solve the linear system to determine these Chern numbers.

\begin{exa} \label{ex:threefold}
Let $I$ be the ideal defined by the $4\times 4$ minors of a $4\times 5$ matrix of general linear forms in $\mathbb C[x_0,x_1,\dots, x_5]$ and let $Z$ be the corresponding threefold in $\p{5}$.
The following table shows the number of solutions resulting from zero-dimensional, 128-bit precision runs in \cite{BHSW} on 5 random elements of $I$ for various choices of degrees. We let $(n_1,n_2,n_3,n_4, n_5)$ denote the degrees of these elements.
\begin{displaymath}
\begin{tabular}{ccccc}
$(n_1,n_2,n_3,n_4,n_5)$ & \vline & non-singular points & \vline & singular points \\ \hline
$(4,4,4,4,4)$ & \vline & 0 & \vline & 1024 \\
$(4,4,4,4,5)$ & \vline & 1 & \vline & 1279 \\
$(4,4,4,5,5)$ & \vline & 6 & \vline & 1594 \\
$(4,4,5,5,5)$ & \vline & 21 & \vline & 1979 \\
\end{tabular}
\end{displaymath}
In each of these computations the zero set consists of the threefold together with a finite set $S$ of points of multiplicity one. The formulas from \cororef{coro:formula} and \thmref{num:criterion} lead to:
\begin{displaymath}
\begin{array}{l}
\deg{c_3}+14\deg{c_2}+61\deg{c_1}+44\deg{Z}=1024 \\
\deg{c_3}+15\deg{c_2}+71\deg{c_1}+65\deg{Z}=1279 \\
\deg{c_3}+16\deg{c_2}+82\deg{c_1}+92\deg{Z}=1594 \\
\deg{c_3}+17\deg{c_2}+94\deg{c_1}+126\deg{Z}=1979.
\end{array}
\end{displaymath}
The unique solution is $\deg{Z}=10$, $\deg{c_1}=0$, $\deg{c_2}=45$, $\deg{c_3}=-46$.
\end{exa}

\begin{exa} \label{ex:threefold2}
Let $V$ be the image of the Segre embedding $i:\mathbb P^1 \times \mathbb P^3 \hookrightarrow \mathbb P^7$. Let $Z$ be the intersection of $V$ with a general hyperplane and let $I$ be the homogeneous ideal of $Z$.
The following table shows the number of solutions resulting from zero-dimensional, 128-bit precision runs in \cite{BHSW} on 6 random elements of $I$ for various choices of the degrees $(n_1,\dots,n_6)$ of these elements.
\begin{displaymath}
\begin{tabular}{ccccc}
$(n_1,n_2,n_3,n_4,n_5, n_6)$ & \vline & non-singular points & \vline & singular points \\ \hline
$(2,2,2,2,2,2)$ & \vline & 0 & \vline & 64 \\
$(2,2,2,2,2,3)$ & \vline & 0 & \vline & 96 \\
$(2,2,2,2,3,3)$ & \vline & 2 & \vline & 142 \\
$(2,2,2,3,3,3)$ & \vline & 10 & \vline & 206 \\
\end{tabular}
\end{displaymath}
In each of these computations the zero set of the 6 random elements consists of the threefold together with a finite set of points of multiplicity one. The formulas from \cororef{coro:formula} and \thmref{num:criterion} lead to:
\begin{displaymath}
\begin{array}{l}
\deg{c_3}+5\deg{c_2}+4 \deg{c_1}-8 \deg{Z}=64 \\
\deg{c_3}+6\deg{c_2}+7 \deg{c_1}-10 \deg{Z}=96 \\
\deg{c_3}+7\deg{c_2}+11 \deg{c_1}-11 \deg{Z}=142 \\
\deg{c_3}+8\deg{c_2}+16 \deg{c_1}-10 \deg{Z}=206.
\end{array}
\end{displaymath}
The unique solution is $\deg{Z}=4$, $\deg{c_1}=10$, $\deg{c_2}=10$, $\deg{c_3}=6$.
\end{exa}

\section{Conclusion}
The results of this paper demonstrate the viability of computing Chern numbers of smooth varieties through numerical homotopy continuation. Homotopy continuation via square systems is a natural venue through which to compute the equivalence of a scheme in the context of an important generalized B\'{e}zout's theorem from intersection theory. It should be noted that the algorithms of this paper could be implemented in a purely symbolic setting as well. The advantages of the numerical approach is that the algorithms work equally well whether the generators are sparse or dense and whether they have rational, algebraic or transcendental coefficients. In addition, meaning can often be attached to the computations in situations where the generators have inaccuracies in their coefficients. Finally, homotopy continuation algorithms parallelize well allowing such algorithms to take full advantage of multi-processor machines. These features suggest the complementary nature of the approach to purely symbolic methods.

In future work, the authors intend to extend the approach to take further advantage of the ideas of intersection theory. In particular, intersection theory for non-square systems leads to a method for computing intersection numbers of Chern classes.

\vskip 20pt

\end{document}